\documentclass{ws-procs9x6}

\newcommand{\ba}{\begin{array}}
\newcommand{\ea}{\end{array}}
\newcommand{\be}{\begin{equation}}
\newcommand{\ee}{\end{equation}}
\newcommand{\Laplacian}{\mathrm{\Delta}}

\newcommand{\pder}[2]{\frac{\partial{#1}}{\partial{#2}}}
\newcommand{\pderr}[2]{\frac{\partial^2 #1}{\partial {#2}^2}}
\newcommand{\R}{\mathbb{R}}

\begin{document}

\title{High order relaxed schemes for nonlinear reaction diffusion problems}

\author{F. Cavalli and M. Semplice$^*$}

\address{Dipartimento di Matematica, Universit\`a di Milano\\
Via Saldini 50, I-20123 Milano, Italy\\
$^*$E-mail: semplice@mat.unimi.it}

\begin{abstract}
Different relaxation approximations to partial differential equations,
including conservation laws, Hamilton-Jacobi equations,
convection-diffusion problems, gas dynamics problems, have been
recently proposed.  The present paper focuses onto diffusive relaxed
schemes for the numerical approximation of nonlinear reaction
diffusion equations.  High order methods are obtained by coupling ENO
and WENO schemes for space discretization with IMEX schemes for time
integration, where the implicit part can be explicitly solved at a
linear cost.  To illustrate the high accuracy and good properties of
the proposed numerical schemes, also in the degenerate case, we
consider various examples in one and two dimensions: the
Fisher-Kolmogoroff equation, the porous-Fisher equation and the porous
medium equation with strong absorption.
\end{abstract}

\keywords{Degenerate reaction diffusion problems, relaxation schemes,
porous-Fisher equation, Fisher-Kolmogoroff equation.}

\section{Introduction}

The main purpose of this work is to approximate solutions of a nonlinear,
possibly degenerate, reaction-diffusion equation of the form
\begin{equation}  \label{eq:degparab}
\frac{\partial u}{\partial t} = D \Laplacian (p(u)) + g(u)
\end{equation}
for $ x\in \Omega \subseteq\mathbb{R}^d,~~d\geq 1,~~ t\geq0$, with
initial condition, $u(x,0)=u_0(x)$ and with suitable boundary
conditions,  
where the function $p$ is a non-decreasing
Lipschitz continuous 
function defined on $\mathbb{R}$. The equation is degenerate if
$p(0)=0$ and in this case the solutions often become non-smooth in
finite time, developing fronts and discontinuities \cite{Aro70}. 
Finally, the coefficient $D$ is a diffusivity
coefficient and the function $g(u)$ is the reaction term.

Equations like \eqref{eq:degparab} are relevant in
describing biological and physical processes and are involved in the modelling of
population growth and dispersal, of waves of concentration of chemical
substances in living organisms, of the motion of viscous fluids (see
e.g. \cite{Murray1}). Also, systems of equations like
\eqref{eq:degparab}, coupled via the reaction terms, are capable of
modelling the cyclic Belousov-Zhabotinskii reactions and the pattern
formation on the wings of butterflies and the coat of mammals (see
\cite{Murray2} and references therein).

This paper is organized as follows. In section
\ref{section:relaxation} we introduce the relaxation approximation of
nonlinear diffusion problems, in section \ref{sec:scheme} we describe
the fully discrete relaxed numerical scheme in the reaction-diffusion
case. In section \ref{sec:tests} we report several numerical tests,
both in one and two space dimensions.

\section{Relaxation approximation of nonlinear diffusion}
\label{section:relaxation}

The schemes proposed in the present work are based on the same
framework of the well-known relaxation approximation of \cite{JX95}
for hyperbolic conservation laws. In the case of the nonlinear
diffusion operator in \eqref{eq:degparab}, an additional variable
$\vec{v}(x,t)\in\mathbb{R}^d$ and a positive parameter $\varepsilon$
are introduced and the following relaxation system is obtained:
\begin{equation}\label{sysrel1}
\left\{
\begin{array}{ll}
\displaystyle
\frac{\partial u}{\partial t} + \mathrm{div}(\vec{v}) = g(u)\\
\\
\displaystyle
\frac{\partial \vec{v}}{\partial t} + \varphi^2 \nabla p(u) = 
    -\frac1\varepsilon \vec{v} + \left(\varphi^2-
    \frac{D}{\varepsilon}\right) \nabla p(u)
\end{array}
\right.
\end{equation}
Now, formally, in the small relaxation limit, $\varepsilon \rightarrow 
0^+$, system \eqref{sysrel1} approximates to leading order
equation \eqref{eq:degparab}.
The parameter $\varphi$ is introduced
in order to have bounded characteristic velocities and to avoid a
singular differential operator as  
$\varepsilon \rightarrow 0^+$.

Finally the non linearity in the convective  term is removed, as in
standard relaxation schemes, introducing another scalar variable
$w(x,t)$ and rewriting the system as:
\begin{equation} \label{3eq}
\left\{
\begin{array}{ll}
\displaystyle
\frac{\partial u}{\partial t} + \mathrm{div}(\vec{v}) = g(u) \\
\\
\displaystyle \frac{\partial \vec{v}}{\partial t} + \varphi^2 \nabla w = 
    -\frac1\varepsilon \vec{v}
    +\left(\varphi^2- \frac{D}{\varepsilon}\right) \nabla w \\
\\
\displaystyle
\frac{\partial w}{\partial t} + \mathrm{div}(\vec{v}) =
    -\frac1\varepsilon (w-p(u))
\end{array} 
\right.
\end{equation}
Formally, as $\varepsilon \rightarrow 0^+$, $w \rightarrow p(u)$, $v
\rightarrow -\nabla p(u)$ and the original equation is recovered. 

In the previous system the parameter $\varepsilon$ has physical
dimensions of time and represents the so-called relaxation
time. Furthermore, $w$ has the same dimensions as $u$, while each
component of $\vec{v}$ has the dimension of $u$ times a velocity;
finally $\varphi$ is a velocity.  The inverse of $\varepsilon$ gives the
rate at which $v$ decays onto $-\nabla p(u)$ in the evolution of the
variable $\vec{v}$ governed by the stiff second equation of
\eqref{3eq}.

Equations \eqref{3eq}, originally introduced in \cite{NP00} for the
purely diffusive case, form a semilinear hyperbolic system with a
stiff source term. The characteristic velocities of the hyperbolic
part are given by $0,\pm\varphi$. The parameter $\varphi$ allows to
move the stiff terms $\frac{D}{\varepsilon}\nabla p(u)$ to the right
hand side, without losing the hyperbolicity of the system.

We point out that degenerate parabolic equations often model physical
situations where free boundaries and discontinuities are relevant: we
expect that schemes for hyperbolic systems will be able to reproduce
faithfully these details of the solution.  One of the main properties
of \eqref{3eq} consists in the semilinearity of the system, that is
all the nonlinear terms are in the (stiff) source terms, while the
differential operator is linear. Hence, the solution of the convective
part requires neither Riemann solvers nor the computation of the
characteristic structure at each time step, since the eigenstructure
of the system is constant in time.  Moreover, the relaxation
approximation does not exploit the form of the nonlinear function $p$
and hence it gives rise to a numerical scheme that, to a large extent,
is independent of it, resulting in a very versatile tool.

We also anticipate here that, in the relaxed case
(i.e. $\varepsilon=0$), the stiff source terms can be integrated
solving a system that is already in triangular form and then it does
not require iterative solvers.

\section{Relaxed numerical schemes}
\label{sec:scheme}

\subsection{Relaxed IMEX schemes}

We extend here the schemes studied in \cite{arxiv0604572}, including
the reaction term. For simplicity, here we describe the
one-dimensional case, the generalization being straightforward.

We observe that system \eqref{3eq} is in the form
\[
 \pder{z}{t} + \pder{f(z)}{x} = g(z) + \frac{1}{\epsilon}h(z)
\]
where $z=(u,v,w)^T$, $f(z)=(v,\varphi^2w,v)^T$, $g(z)=(g(u),0,0)^T$
and $h(z)= (0,-v+(\varepsilon\varphi^2-D)w_x,p(u)-w)^T$. When
$\varepsilon$ is small, the presence of both 
non-stiff and stiff terms, suggests the use of IMEX schemes
\cite{ARS97,KC03,PR05}. In the problems considered here, the reaction
term $g(u)$ is not stiff and hence we treat it with the explicit
portion of the IMEX scheme.

First we consider a semi-discrete form for time integration.
Let's assume for simplicity a uniform time step $\Delta t$ and
denote with $z^n$ the numerical approximation of the variable $z$ at
time $t_n=n\Delta t$, 
for $n=0,1,\ldots$ In our case a $\nu$-stages IMEX scheme reads
\be \label{eq:IMEX}
  z^{n+1} = z^n 
     - \Delta t \sum_{i=1}^{\nu}
           \tilde{b}_i  \left[\pder{f}{x}(z^{(i)}) +g(z^{(i)})\right]
     + \frac{\Delta t}{\varepsilon} \sum_{i=1}^{\nu} b_i h(z^{(i)})
\ee
where the stage values are computed as
\be \label{eq:stage:i}
z^{(i)} = B^{(i)}
      + \frac{\Delta t}{\varepsilon} a_{i,i} h(z^{(i)})
\ee
for
\be \label{eq:stage:expl}
B^{(i)} = z^n 
        -\Delta t \sum_{k=1}^{i-1}\tilde{a}_{i,k} \left[\pder{f}{x}(z^{(k)}) +g(z^{(k)})\right]
        + \frac{\Delta t}{\varepsilon} \sum_{k=1}^{i-1} a_{i,k} h(z^{(k)})
\ee
Here $({a}_{ik},{b}_i)$ and  $(\tilde{a}_{ik},\tilde{b}_i)$ are a pair
of Butcher's tableaux \cite{HW1} of, respectively, a diagonally implicit and an
explicit Runge-Kutta schemes.

In this work we use the so-called relaxed schemes, that are obtained
by formally letting $\varepsilon\rightarrow0$ in \eqref{eq:IMEX}. For
these the computation of the first stage, that is \eqref{eq:stage:i}
with $i=1$,  
\be \label{eq:stage:1} 
 \left[\begin{array}{l} 
    {u^{(1)}}\\
    {v^{(1)}}\\
    {w^{(1)}}
 \end{array}\right]
 = 
 \left[\begin{array}{l} u^{n}\\v^{n}\\w^{n}\end{array}\right]
 +
 \frac{\Delta t}{\varepsilon} a_{1,1} 
   h\left(\left[\begin{array}{l}u^{(1)}\\v^{(1)}\\w^{(1)}\end{array}\right]\right)
\ee
implies that $h(z^{(2)})=0$, which is equivalent to
\[ {u^{(1)}=u^n} 
   \qquad {w^{(1)}=p(u^{(1)})} 
   \qquad v^{(1)}=\pder{w^{(1)}}{x} ,\]
Now the second stage, $i=2$, reads 
\[z^{(2)} = z^n 
    - \Delta{t}\tilde{a}_{2,1}{\left[\pder{f}{x}(z^{(1)}) + g(z^{(1)})\right]} 
    + \frac{\Delta t}{\varepsilon} a_{2,1}
           \underbrace{{h(z^{(1)})}}_{{\equiv 0}}
    + \frac{\Delta t}{\varepsilon} a_{2,2} {h(z^{(2)})}.
\]
Hence the last two components of $z^{(2)}$ are determined by the stiff
terms of the above expression, namely $h(z^{(1)})=0$. On the other
hand, due the form of
$h(z)$, there are no stiff terms in the equation for the first
component $u^{(2)}$, which  is then determined by a balance law.

Summarizing, the relaxed scheme yields an alternation of
{\em relaxation steps}
\[
  h(z^{(i)})=0 \qquad 
   \text{ i.e. }\, 
     w^{(i)}=p(u^{(i)}),
     v^{(i)}=\pder{w^{(i)}}{x}
\]
and {\em transport steps} where we advance for time $\tilde{a}_{i,k}\Delta{t}$
\be \label{eq:transportRK}
 \pder{z}{t} + \pder{f(z)}{x} = g(z) 
\ee
with initial data $z=z^{(i)}$, retain only the first component and
assign it to $u^{(i+1)}$.
Finally the value of $u^{n+1}$ is computed as \(u^n+\sum \tilde{b}_i
u^{(i)}\).

\subsection{Spatial reconstructions}

In order to have a fully discrete scheme,
we still need to specify
the space discretization. We will use discretizations based on finite
differences, in order to avoid cell coupling due to the source terms.

Recall that the IMEX technique reduces the integration to a cascade of
relaxation and transport steps. The former are the implicit parts of
\eqref{eq:stage:1} and \eqref{eq:stage:i}, while the transport steps
appear in the evaluation of the explicit terms $B^{(i)}$ in
\eqref{eq:stage:expl}. Since \eqref{eq:stage:1} and \eqref{eq:stage:i}
involve only local operations, the main task of the space discretization
is the evaluation of $\partial_x{f}$, where we will exploit the
linearity of $f$ in its arguments.

In the one-dimensional case, let us consider a uniform grid on $[a,b]\subset\mathbb{R}$,
\(x_j=a-\frac{h}2+jh\) for $j=1,\ldots,m$, where $h=(b-a)/m$ is the
grid spacing and $m$ the number of cells. 
We denote with $z^j_n$ the value of the quantity $z$ at time $t^n$ at
$x_j$, the centre of the $j^{\mbox{\scriptsize th}}$ computational cell. 
The fully discrete scheme
may be written as 
\[
z_j^{n+1} = z_j^n - \Delta t \sum_{i=1}^{\nu}
     \tilde{b}_i  \left(F^{(i)}_{j+1/2} - F^{(i)}_{j-1/2}\right)
    + \frac{\Delta t}{\varepsilon} \sum_{i=1}^{\nu} b_i
     g(z_j^{(i)}),  
\]
where $F^{(i)}_{j+1/2}$ are the numerical fluxes, which are the only
item that we still need to specify. It is necessary to
write the scheme in conservation form and thus, following
\cite{OS89}, we introduce the function $\hat{F}$ such that
\[
  f(z(x,t))=\frac1h \int_{x-h/2}^{x+h/2} \hat{F}(s,t) \mathrm{d}s
\]
and hence
\[
  \frac{\partial f}{\partial x}(z(x_j,t)) = \frac1h
  \left(\hat{F}(x_{j+1/2},t)-\hat{F}(x_{j-1/2},t)\right).
\]
The numerical flux function $F_{j+1/2}$ approximates
$\hat{F}(x_{j+1/2})$. 

In order to compute the numerical fluxes, for each stage value, we
reconstruct boundary extrapolated data $z^{(i)\pm}_{j+1/2}$ with a
non-oscillatory interpolation method, starting from the point values
$z^{(i)}_j$ of the 
variables at the centre of the cells. Next we apply a monotone
numerical flux to these boundary extrapolated data.

To minimize numerical viscosity we choose the Godunov flux, which, in
the present case of a linear system of equations, reduces to the upwind
flux. In order to select the upwind direction we write the linear
system with constant coefficients
\eqref{eq:transportRK} in
characteristic form. The characteristic variables relative to the
eigenvalues $\varphi,-\varphi,0$  are respectively
\[
      U=\frac{v+\varphi w}{2\varphi}     \qquad      
      V=\frac{\varphi{}w-v}{2\varphi}  \qquad   
      W=u-w.
\]
Note that $u=U+V+W$. Therefore the numerical flux in characteristic
variables is
\(F_{j+1/2}=(\varphi U^-_{j+1/2},-\varphi V^+_{j+1/2},0)\).

The accuracy of the scheme depends on the accuracy of the
reconstruction of the boundary extrapolated data. For a first order
scheme we use a piecewise constant reconstruction such that
\(U^-_{j+1/2}=U_j\) and \(V^+_{j+1/2}=V_{j+1}\). For higher order
schemes, we use ENO or WENO reconstructions of appropriate accuracy
\cite{OS89}.
 
Since the transport steps need to be applied only to
$u^{(i)}$, we have 
\[
 u_j^{(i)}
 = u_j^n 
   -\lambda \sum_{k=1}^{i-1}\tilde{a}_{i,k}
   \left[
   \varphi \left( U^{(k)-}_{j+1/2}-U^{(k)-}_{j-1/2}
   -V^{(k)+}_{j+1/2}+V^{(k)+}_{j-1/2}\right)
   +\mathrm{\Delta t} g(u^{(k)}_j)
   \right]
\]
Finally, taking the last stage value and going back to conservative
variables, 
\[
\ba{ll} 
u_j^{n+1}=u_j^n 
    -\frac\lambda2
     \sum_{i=1}^\nu \tilde{b}_i 
         &\left([v^{(i)-}_{j+1/2}+v^{(i)+}_{j+1/2}-(v^{(i)-}_{j-1/2}+v^{(i)+}_{j-1/2})]\right.\\
              &\left.+\varphi[w^{(i)-}_{j+1/2}-w^{(i)+}_{j+1/2}-(w^{(i)-}_{j-1/2}-w^{(i)+}_{j-1/2})]
         \right)
\ea
\]

We wish to emphasize that the scheme reduces to the time
advancement of the single variable $u$. Although the scheme is based
on a system of three equations, the construction is used only to select
the correct upwinding for the fluxes of the relaxed scheme and the
computational cost of each time step is not affected.

\subsection{Numerical scheme}

Employing a Runge-Kutta IMEX scheme of
order $p$, can give an integration procedure for \eqref{eq:degparab}
which is of order up to $2p$ with respect to $h$ (see
\cite{arxiv0604572}),
because the CFL restrictions on the time-step are of parabolic type
($\mathrm{\Delta{}t}\leq{}Ch^2$). We observed that this
theoretical convergence rate can be achieved in practice, with careful
choice of the approximations of the spatial derivatives.

Summarizing, the relaxed schemes that we propose for the numerical
integration of \eqref{eq:degparab} consists of the following
steps. 

For each Runge-Kutta stage we need to compute the variables
$v^{(i)}$ and $w^{(i)}$ (relaxation steps). The computation of
$v^{(i)}$ requires the approximation of a spatial gradient operator,
for which we choose a central finite difference operator of order at
least $2p$. Then we need to
solve the transport equation by diagonalizing the linear system
\eqref{eq:transportRK}, reconstructing the characteristic variables
$U^{(i)}$ and $V^{(i)}$ at cell boundaries and computing the
fluxes. Again, we must choose a spatial reconstructions of order at
least $2p$ and, to avoid spurious
oscillations, we employ ENO or WENO non-oscillatory procedures. These
procedures compare, for each cell, the reconstructions obtained using
different stencils and choose the least oscillatory one (ENO) or
compute a weighted linear combination of all of them (WENO). For
details, see \cite{OS89}.

\section{Numerical results}
\label{sec:tests}

\subsection{Travelling waves tests}

\begin{figure}
\psfig{file=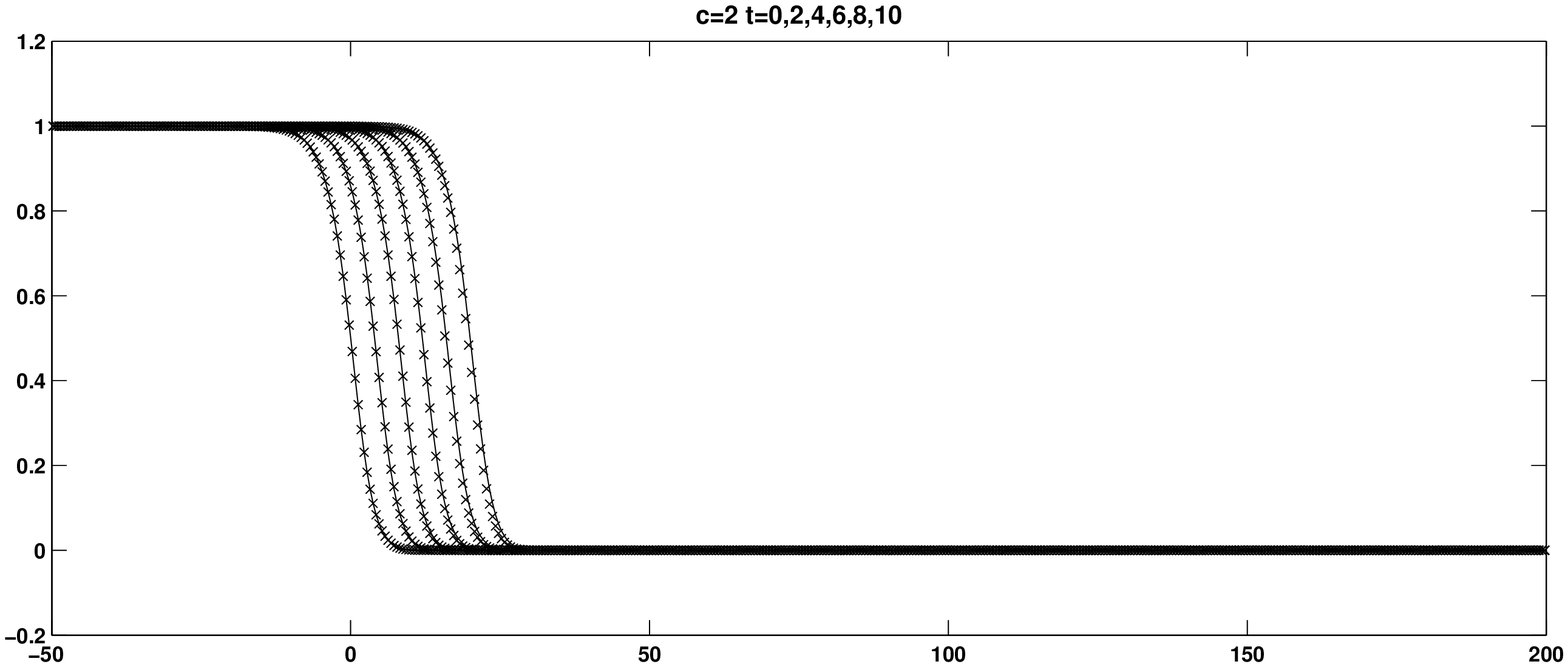,width=\textwidth}

\psfig{file=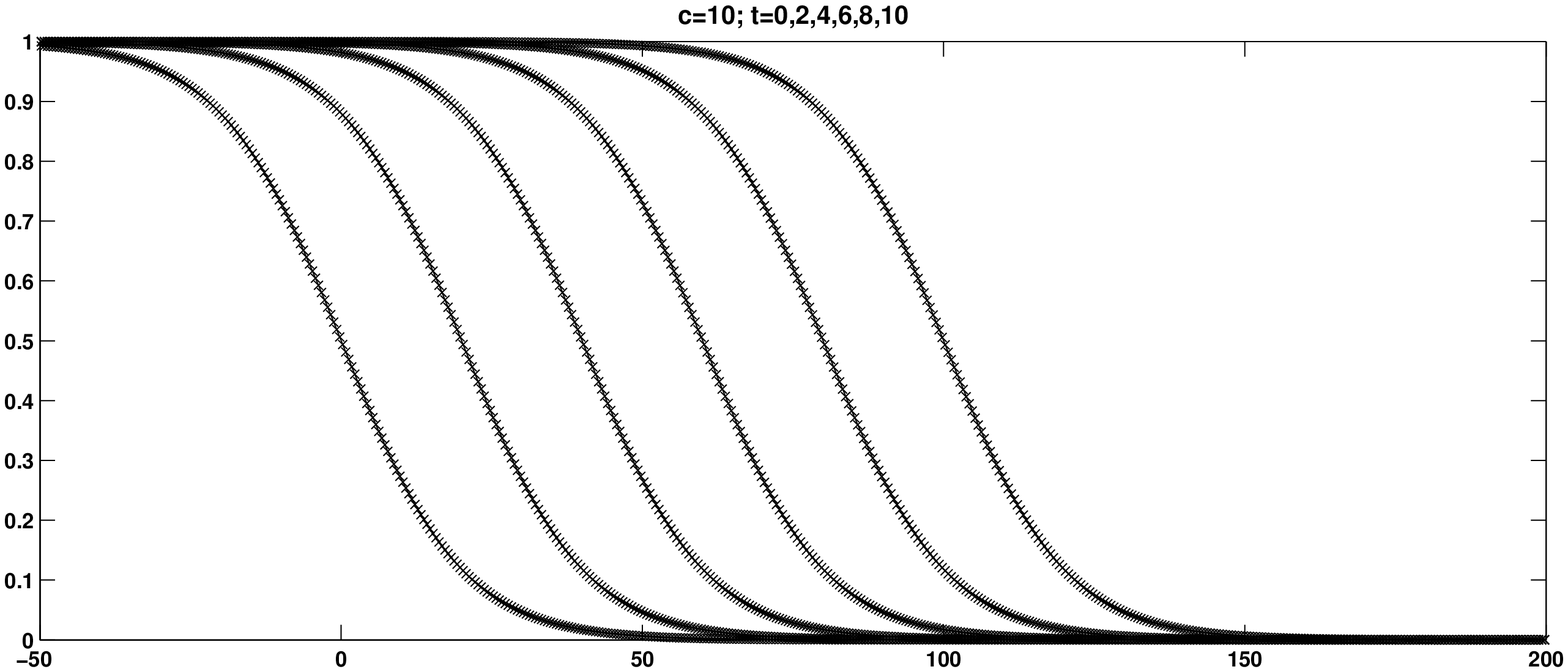,width=\textwidth}
\caption{Travelling waves with $c=2$ (top) and $c=10$ (bottom). The
  crosses represent the numerical solution at $t=0,2,4,6,8,10$ and the
  solid lines are the corresponding asymptotic expansions
  \eqref{eq:fk:tw}.}
\label{fig:tw}
\end{figure}

As a first test, we consider the one-dimensional Fisher-Kolmogoroff equation, namely
\begin{equation} \label{eq:fk}
 \pder{u}{t} = ku(1-u) + D \pderr{u}{x} 
\end{equation}
for $x\in\R$ and $t>0$, which was initially proposed for the modelling of spatial
spread of populations \cite{Murray1}.

The two uniform steady states of \eqref{eq:fk} are $u_0(x)=0$ and
$u_1(x)=1$. A careful analysis of the state space, reveals the
existence of travelling wave solutions linking these two
states, i.e. of solutions of the form $u(x,t)=U(x-ct)$ such that
$U(-\infty)=1$ and $U(\infty)=0$. \cite{Murray1} reports the following
asymptotic form for such solutions:
\begin{equation} \label{eq:fk:tw}
 U(z)=
 \frac1{(1+e^{z/c})} +
 \frac1{c^2}\frac{e^{z/c}}{(1+e^{z/c})^2}\log\frac{4e^{z/c}}{(1+e^{z/c})^2} +
 O\left(\frac1{c^4}\right),
\end{equation}
which is valid when the speed verifies
\[ c\geq 2\sqrt{kD}. \]
Moreover the slope of the inflection point ($z=0$) is related to the
speed $c$ of the wave by the relation
\begin{equation} \label{eq:fk:sc}
U'(0) = -\frac1{4c} + O\left(\frac1{c^5}\right), \mbox{for } c\geq2.
\end{equation}
Travelling wave solutions are very important in the modelling of
populations, since they represents phenomena like the invasion of a
territory by a new species, the expansion of epidemics, etc.

\begin{figure}
\begin{center}
\psfig{file=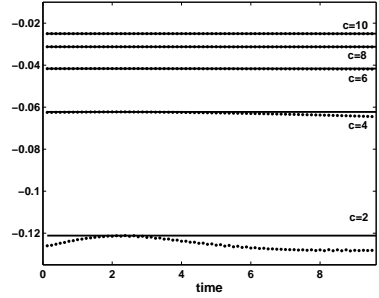,width=.5\textwidth}
\end{center}
\caption{Comparison of the maximum slope of the numerical solutions
  (dots) with the value predicted by the asymptotic expansion
  \eqref{eq:fk:sc}. Time is on the horizontal axis, the values of
  $c$ are printed on the graph close to each line.} 
\label{fig:slopespeed}
\end{figure}

We tested the ability of our integration scheme to reproduce these
results. We chose $k=D=1$ in \eqref{eq:fk}, set up initial datum
$u(x,0)=U(x)$ for $c=2,4,\ldots,10$ on $x\in[-50,200]$ and evolved it
with homogeneous Neumann boundary conditions until $T=10$. We present
the numerical solutions obtained with ENO spatial reconstruction of
order $6$ and third order Runge-Kutta time integration: they are
compared with the asymptotic ones in Figure \ref{fig:tw}.  Moreover we
tested the validity of the expansion \eqref{eq:fk:sc} plotting the
maximum of the numerical gradient of the solution, together with the
values predicted by the asymptotic analysis. From Figure
\ref{fig:slopespeed} we clearly see a very accurate agreement, except
for the limit case $c=2$, which is however at the boundary of the
validity range of \eqref{eq:fk:sc}.

\begin{figure}
\begin{center}
\psfig{file=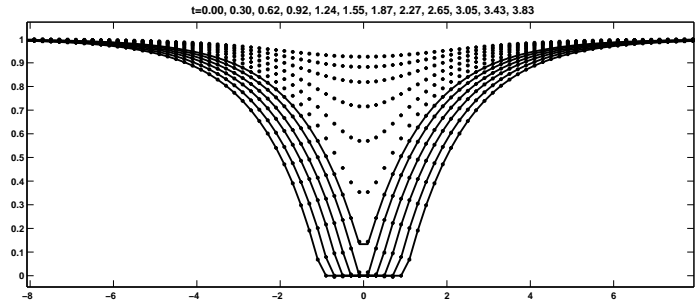,width=\textwidth}
\end{center}
\caption{Comparison of the numerical solutions (dots) of \eqref{eq:genfk} with
  the asymptotic expansion (solid line). Here $p=1$ and
$q=m$. The latter is printed only in
  the time range where it is valid.} 
\label{fig:wit}
\end{figure}

As a second test we consider the following generalization of the
Fisher-Kolmogoroff equation \eqref{eq:fk}
\begin{equation} \label{eq:genfk}
 \pder{u}{t} = u^p(1-u^q) +  \pder{}{x}\left(u^m\pder{u}{x}\right) 
\end{equation}
The existence of travelling waves can be proved for a wide range of
parameters $p,q,m$ \cite{Murray1}. The paper
\cite{Wit95} gives an expression of such waves for the case $p=1$ and
$q=m$. Moreover \cite{Wit95} finds an asymptotic expansion for
the case of two merging travelling waves, which is valid for a finite
time after the first contact of the two fronts.

Figure \ref{fig:wit} shows both the numerical
solutions obtained with our method and the asymptotic expansion given
in \cite{Wit95}, in the case $q=m=1$. 
The initial data represent two travelling fronts, initially at $x_0=-1$
and $x_1=1$, moving in opposite
directions. They first meet at $t^*=1.41$ and then the two waves
merge. The asymptotic expansion given in \cite{Wit95} is valid for a
small time interval after $t^*$ and it is shown in Figure
\ref{fig:wit} with solid lines. The numerical solutions (in dots) are
shown until the steady state is reached.

\subsection{Two dimensional tests}

\begin{figure}
\begin{center}
\psfig{file=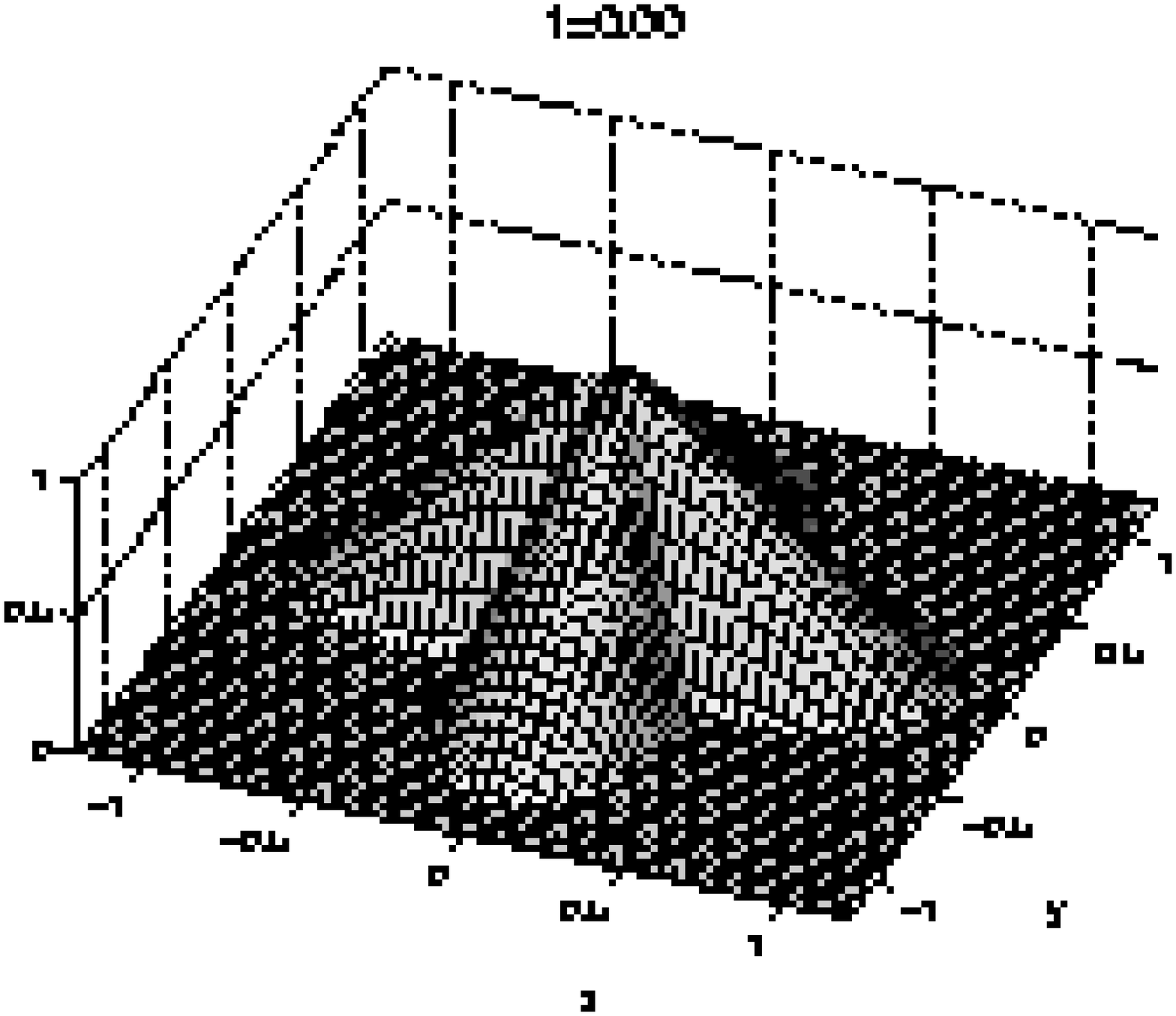,width=.45\textwidth}
\psfig{file=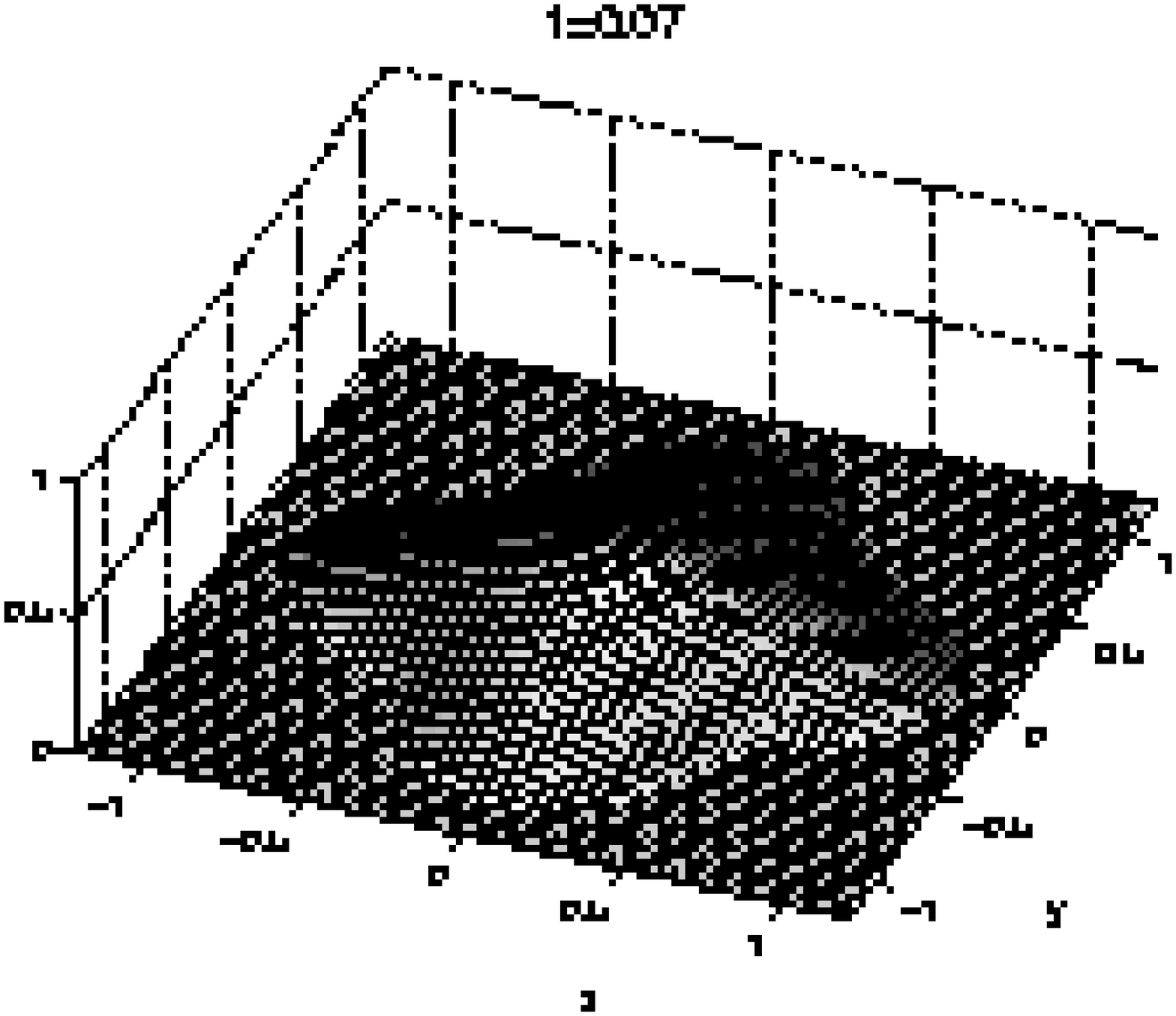,width=.45\textwidth}

\psfig{file=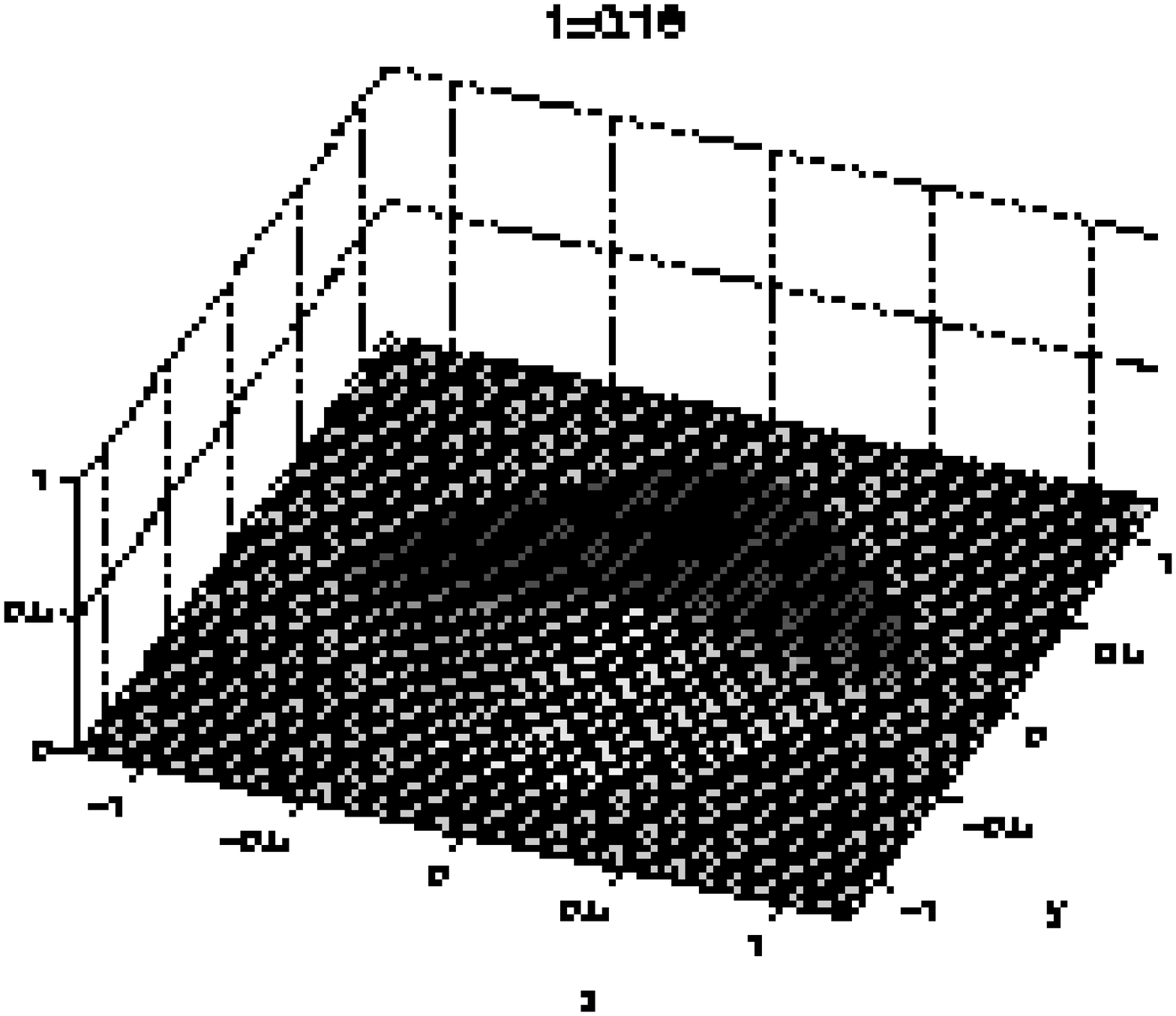,width=.45\textwidth}
\psfig{file=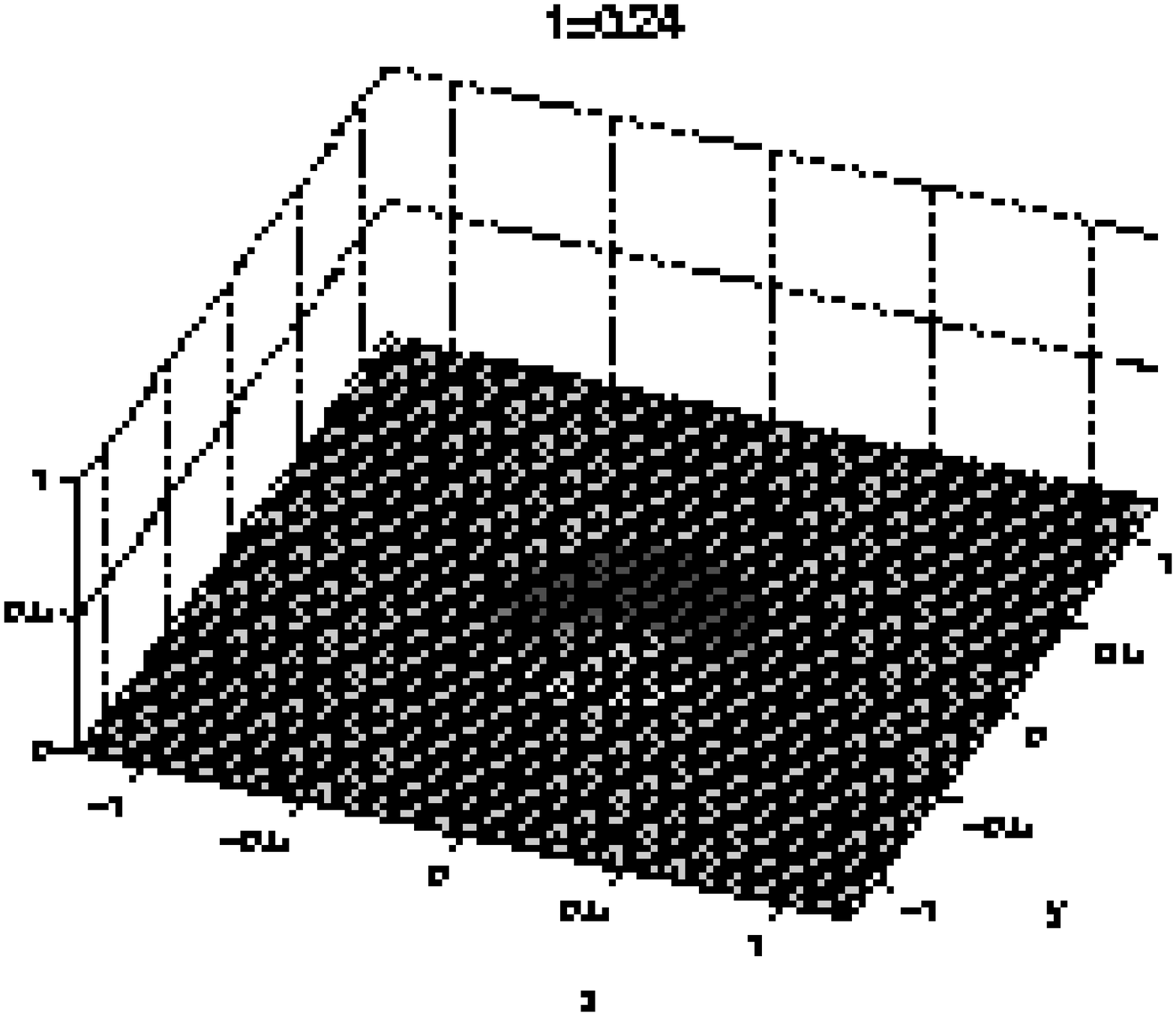,width=.45\textwidth}
\end{center}

\caption{The numerical solution of \eqref{eq:mik} with initial data
  \eqref{eq:croce} at different times, for
  $(x,y)\in[-1.2,1.2]^2$.}
\label{fig:testcroce}
\end{figure}

\begin{figure}
\psfig{file=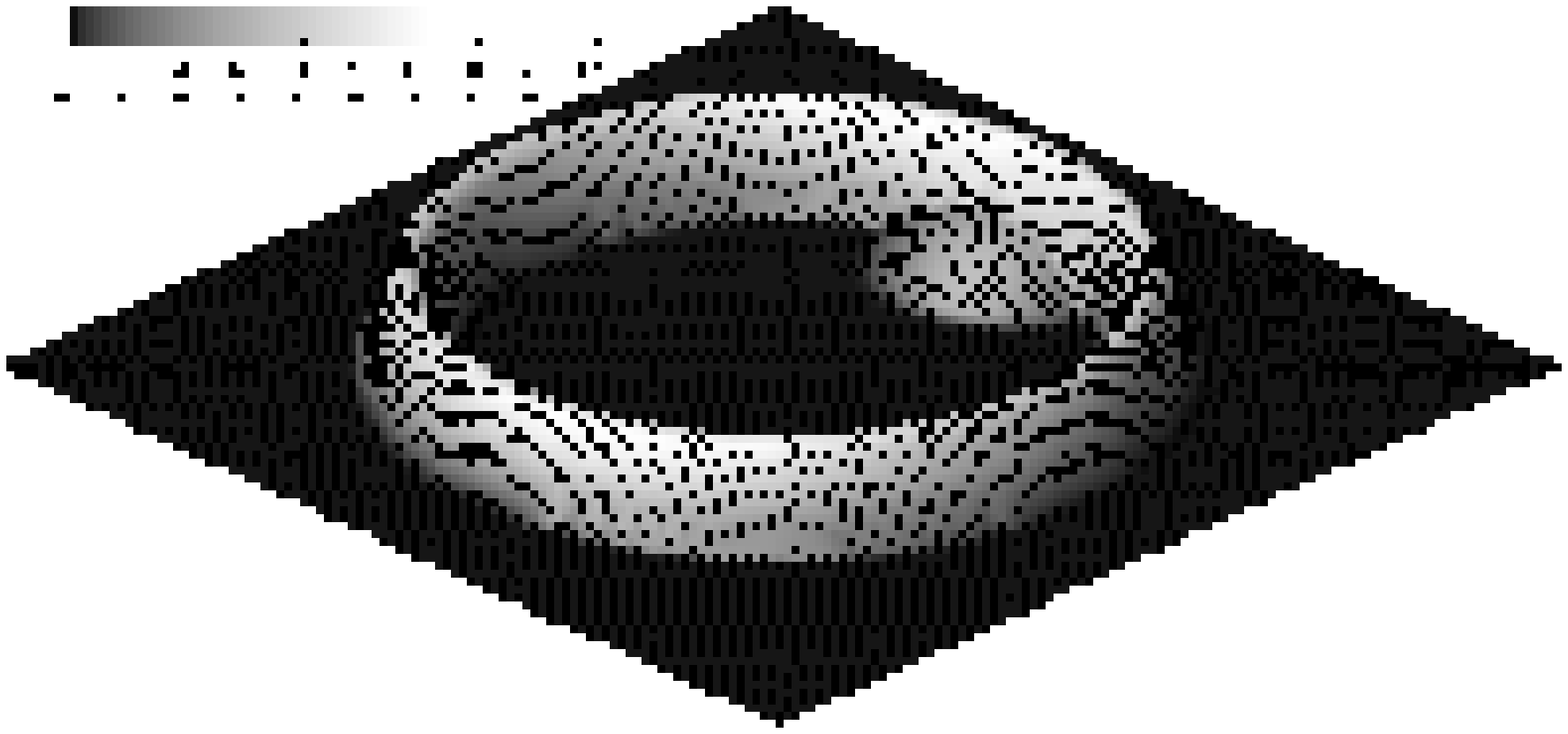,angle=-90,width=.95\textwidth,}

\psfig{file=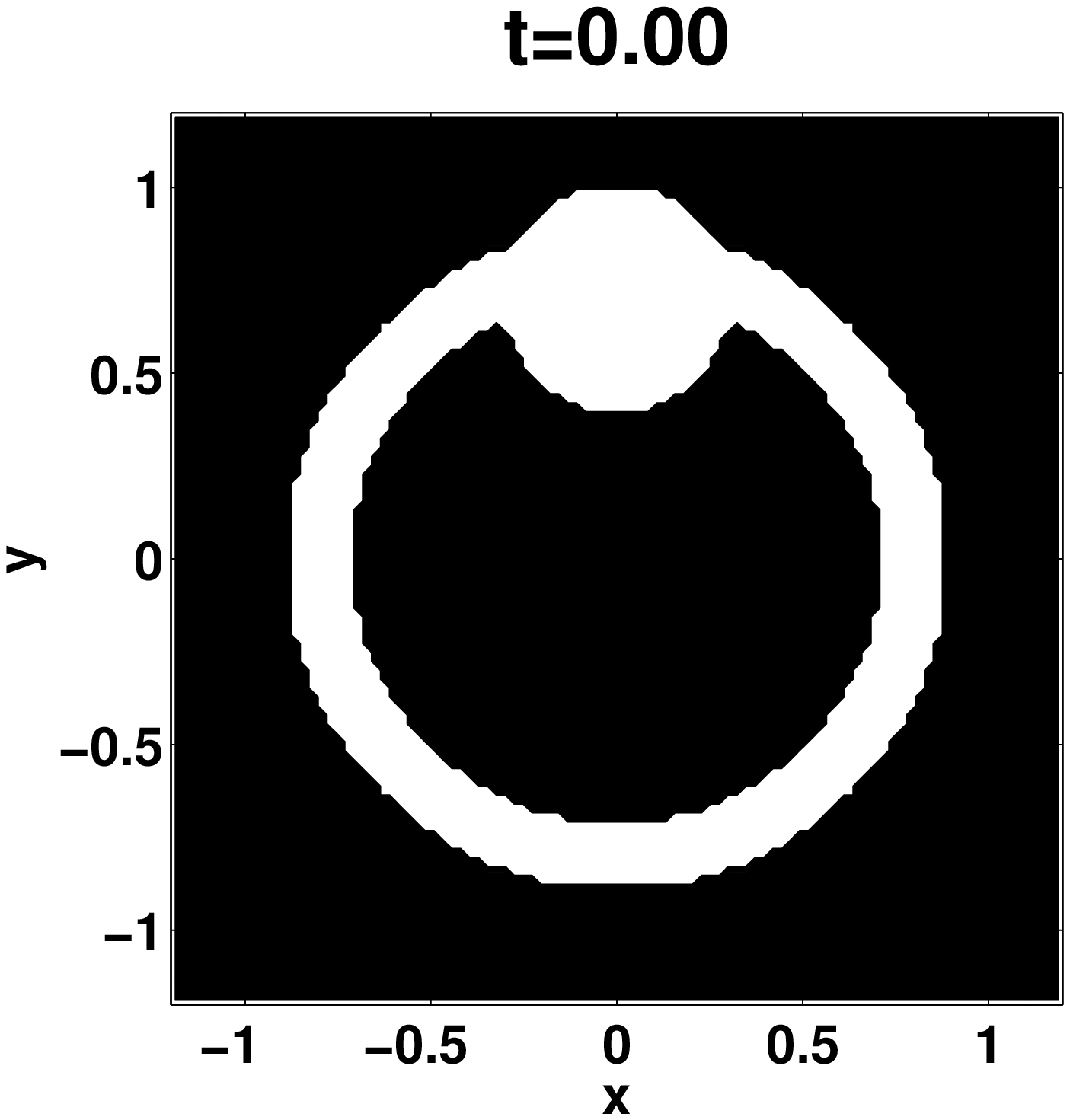,width=.24\textwidth}
\psfig{file=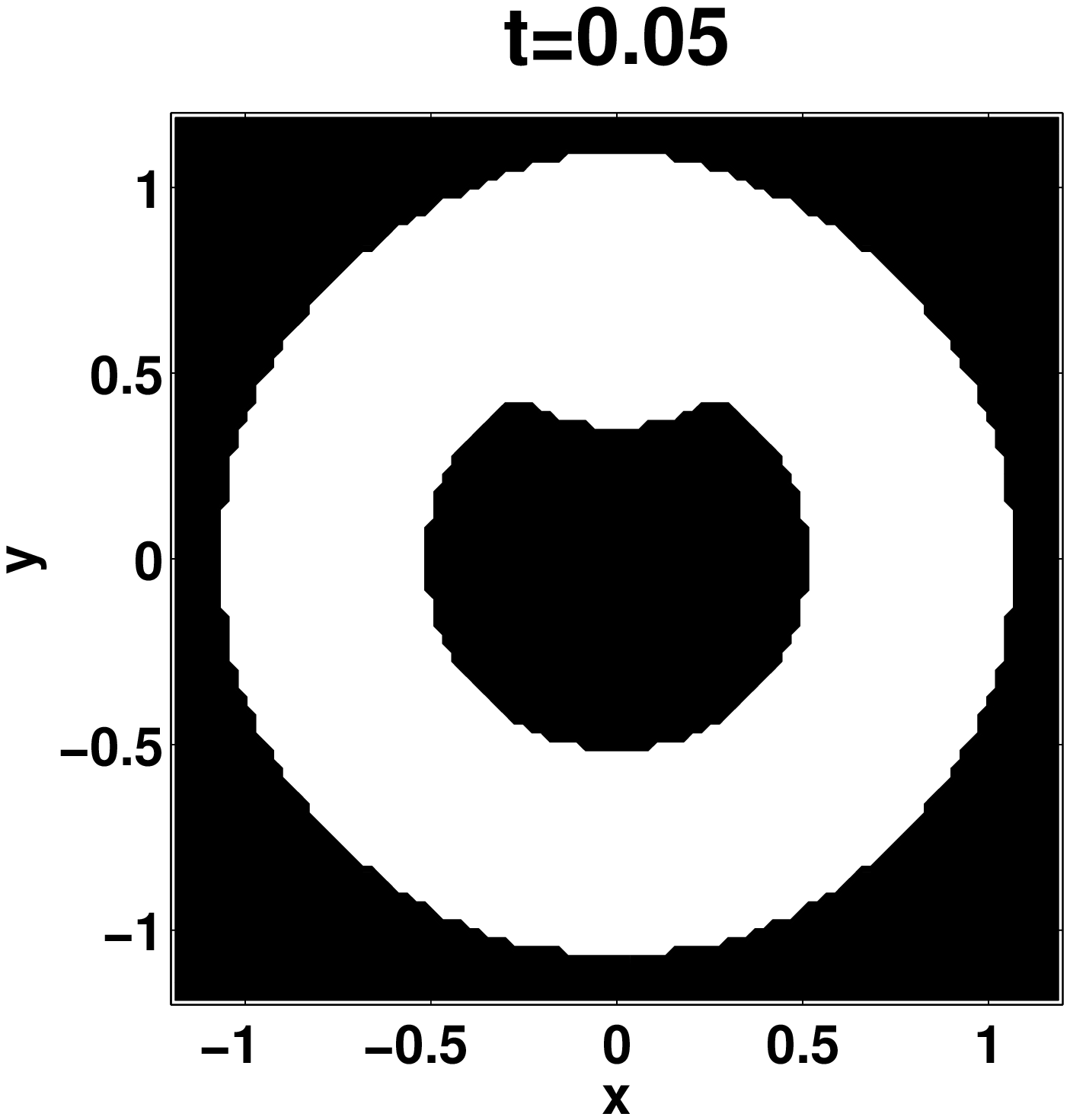,width=.24\textwidth}
\psfig{file=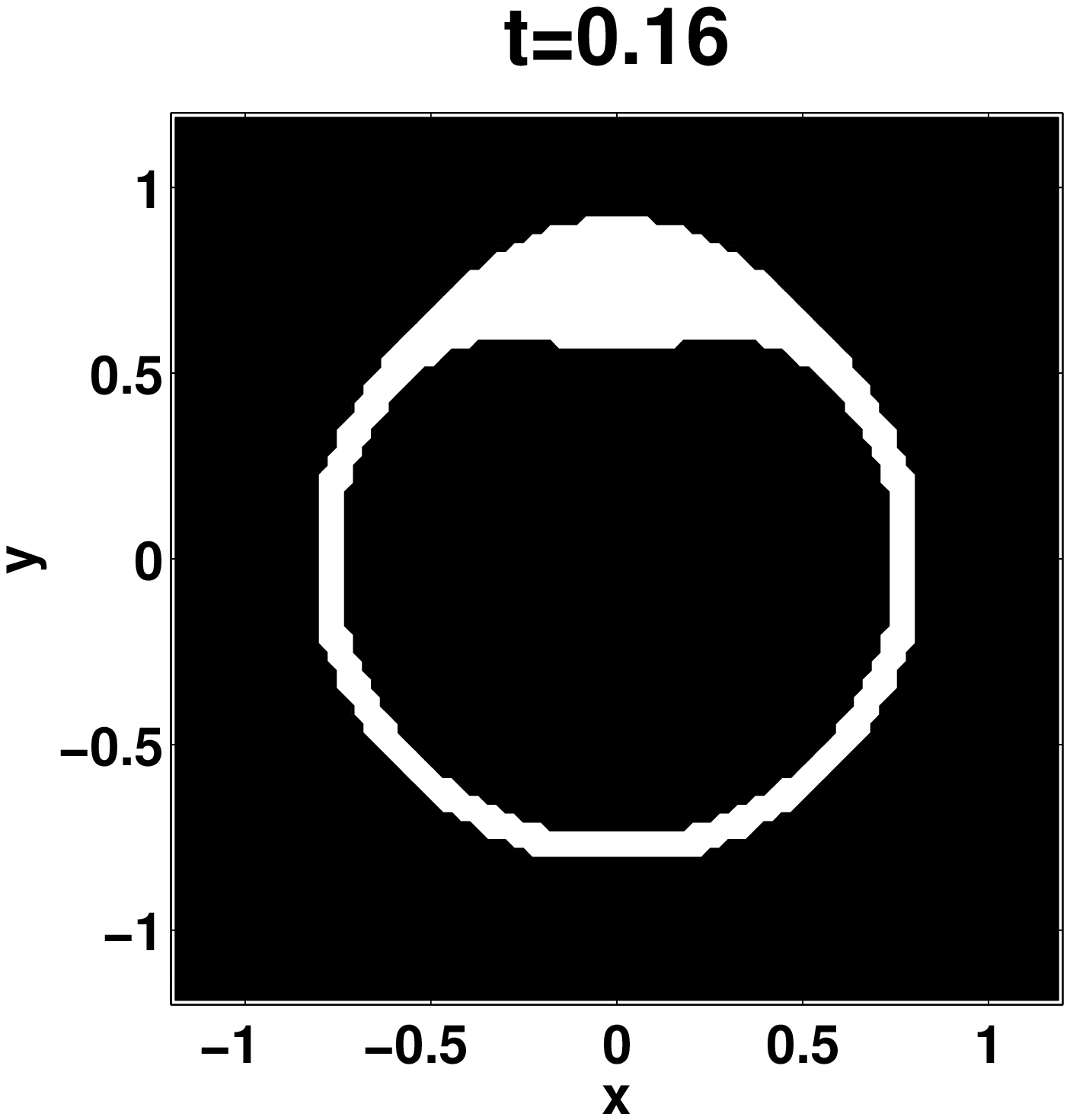,width=.24\textwidth}
\psfig{file=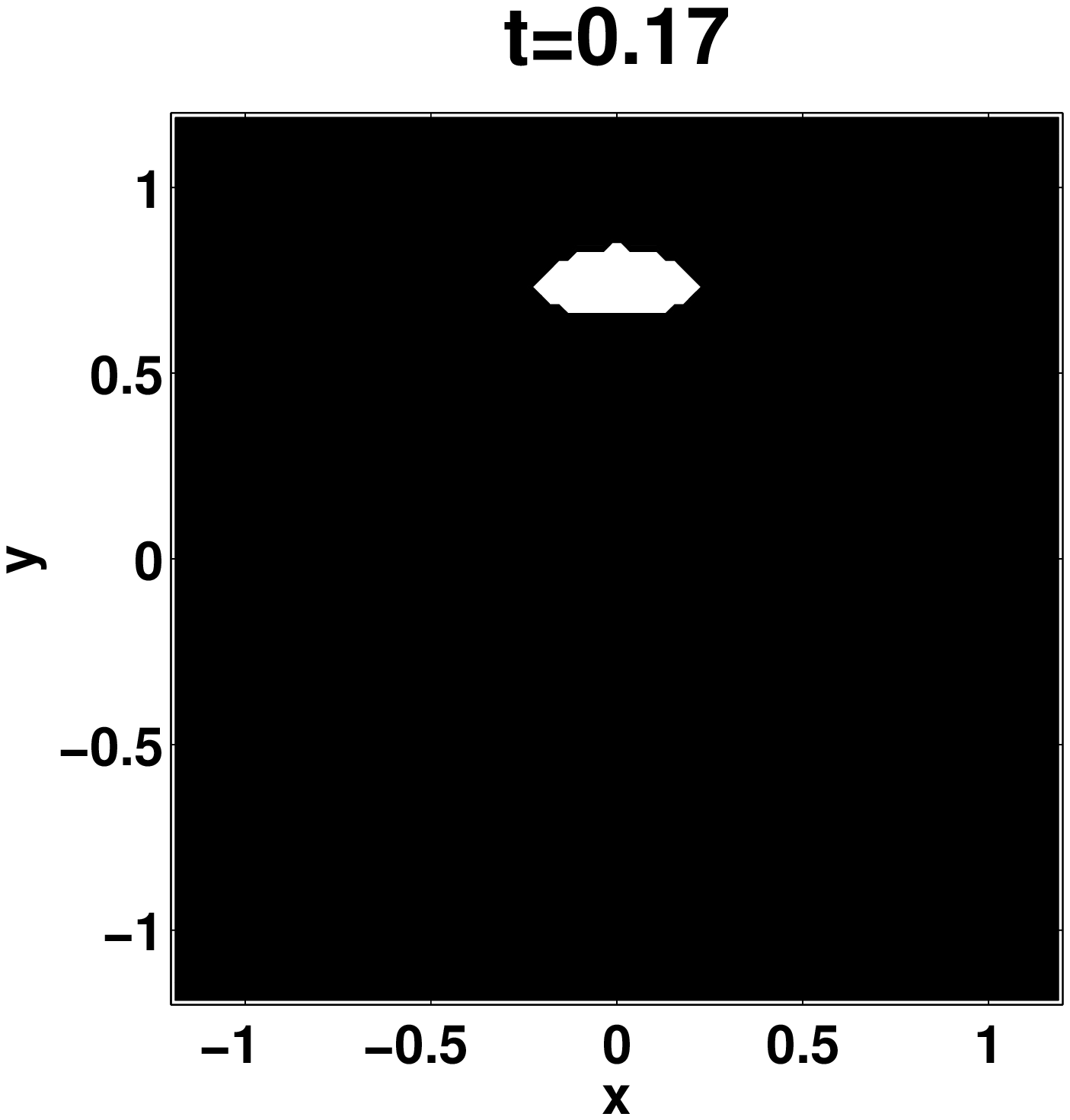,width=.24\textwidth}

\caption{Initial data (top) and supports (bottom) of the numerical solution of
  \eqref{eq:mik} at different times.}
\label{fig:testanellobolla}
\end{figure}

We tested the multi-dimensional version of our integration scheme on
the equation
\begin{equation} \label{eq:mik}
 \pder{u}{t} = \mathrm{\Delta}(u^m) -cu^p 
\end{equation}
for $x\in[-2,2]^2$ and $t\geq0$. The above equation presents
interesting finite-time extinction phenomena, reported in
\cite{Mik95}. We first take initial data
\begin{equation} \label{eq:croce}
 u(x,y,0) = 
   \left\{
   \begin{array}{ll}
   1 & ,x=0,y=0\\
   \left[1-\frac{(x^2+y^2)^2}{\sqrt{x^6+y^6}} \right]_+ &
   \mbox{otherwise}
   \end{array}
   \right.
\end{equation}
and evolve it until extinction with equation \eqref{eq:mik} for $c=5$,
$p=0.5$ and $m=2$. The results are presented in Figure
\ref{fig:testcroce}. 

Finally we tested the persistence of asymmetry in the initial datum
along the evolution. We took $u(x,y,0)$ as a radially symmetric
function with a small perturbation, see
Figure \ref{fig:testanellobolla}, and evolved it with the same equation
and parameters as before, until extinction. Figure
\ref{fig:testanellobolla} shows clearly that the initial perturbation
of the radial symmetry is maintained until the solution vanishes.

\section{Conclusions}

We have proposed and analyzed relaxed schemes for nonlinear degenerate
reaction diffusion equations.
By using suitable discretization in space and time, namely ENO/WENO
non-oscillatory reconstructions for numerical fluxes and IMEX
Runge-Kutta schemes for time integration, we have obtained a class of
high order schemes. The theoretical convergence
analysis for the semidiscrete scheme and the
stability for the fully discrete schemes have been studied by us, for
the case of nonlinear diffusion, in
\cite{arxiv0604572}.

Here we tested these schemes on travelling waves solutions and in
cases where the solution vanishes in finite time. In all cases we
observed a very good agreement with known properties of the exact
solutions.

\bibliographystyle{ws-procs9x6}
\bibliography{CavalliSemplice}

\end{document}